# General hypergeometric distribution: A basic statistical distribution for the number of overlapped elements in multiple subsets drawn from a finite population

Xing-gang Mao[1,3*] and Xiao-yan Xue[2,3]

[1]Department of Neurosurgery, Xijing Hospital, the Fourth Military Medical University, Xi'an, No. 17 Changle West Road, Xi'an, Shaanxi Province, China;

[2]Department of Pharmacology, School of Pharmacy, the Fourth Military Medical University, No. 17 Changle West Road, Xi'an, Shaanxi Province, China

[3]These authors contributed equally to this work.

*Correspondence:xgmao@hotmail.com (X.G.M.)

Short title: General hypergeometric distribution

**Abstract**

General hypergeometric distribution (GHGD) describes the following distribution: from a finite space containing N elements, select T subsets with each subset contains M[i] (T-1 $\geq$ i $\geq$ 0) elements, what is the probability that exactly x elements are overlapped exactly t times or at least t times ($X_{LO=t}$ or $X_{LO\geq t}$, T $\geq$ t $\geq$ 0, here LO is level of overlap)? The classical hypergeometric distribution (HGD) describes the situation of two subsets, while the general situation has not been resolved, despite the overlapped elements has been visualized with the Venn diagram method for about 140 years. GHGD described not only the distribution of $X_{LO=t}$ or $X_{LO\geq t}$ that are overlapped in all of the subsets ($X_{LO=T}$), but also the $X_{LO=t}$ or $X_{LO\geq t}$ that are overlapped in a portion of the subsets (LO = t or LO $\geq$ t, T $\geq$ t $\geq$ 0). Here, we developed algorithms to calculate the GHGD and discovered graceful formulas of the essential statistics for the GHGD, including mathematical expectation, variance, and high order moments. In addition, statistical theory to infer a statistically reliable gene set from multiple datasets based on these formulas was established by applying Chebyshev's inequalities.

## 1 Introduction

Statistically dealing with finite populations are commonly difficult tasks (Casella and Berger, 2002). Suppose selecting *K* balls at random from a space of totally *N* balls containing *M* red and *N-M* green balls, the classical hypergeometric distribution (HGD) describes the probability that exactly *x* of the selected balls are red $p(X = k \mid N, M, K)$. If we consider the *M* red balls and the selected *K* balls as two subsets of the space, then the selected red balls are the intersection of the two subsets. Therefore, the classical HGD described the number of overlapped elements (NOE) in two subsets (containing *M* and *K* elements, respectively) selected from a finite population space (containing *N* elements). However, if the selected subsets are more than *2*, there is yet no statistical distribution can describe the NOE in these subsets. Suppose there are totally *T* (T≥2) subsets, then we define these kinds of distribution as *general hypergeometric distribution* (GHGD): from a finite space containing *N* elements, select *T* subsets with each subset contains M[*i*] elements (*T-1* ≥ *i* ≥ 0), what is the probability that exactly *x* elements are overlapped exactly *t* times or at least *t* times ($X_{LO=t}$ or $X_{LO \geqslant t}$, T ≥ t ≥ 0, here *LO* is **level of overlap**)? This issue can also be described as the problem of group allocation of particles, for which the number of empty cells (corresponding to the number of non-overlapped elements, $X_{LO=0}$) has been studied (Vatutin and Mikhajlov, 1982).

The GHGD is important because it describes the statistical behavior of random subsets selections, and is commonly visualized by the Venn diagram method. Because set is a fundamental concept in modern science, the statistical properties of families of subsets would be a basic theory for many areas such as set theory, combinatorics, probability, graph theory, group theory, algebraic topology (where subsets are subcomplex, and overlapped elements are common vertices), et al. In practice, there is also a huge potential demand of such a theory along with the accelerated accumulation of big data. A remarkable area dealing with this situation more and more widely is in biological research, especially in the context of rapid progress of high throughput data (genome, proteome, et al) acquisition technology (Table S1 listed a small part of works using the methods in recent years). The elements with a certain degree of overlap in multiple sets sampled from a finite population space (normally genome, proteome, et al) were often considered to be worth for further study[2, 3], which is used in a wide range of situations [4-8], and played essential roles for some studies [9-11]. However, although the Venn diagram method has been used for about 140 years, the statistical distribution has not been established yet for the NOE, mainly because of the difficulty of the problem, as partially reflected by the challenge to develop tools to plot the Venn diagram [12-15] even for very small amount of sets (such as *4 - 6* sets). Because of lacking a statistical distribution theory for GHGD, a few studies used classical HGD for two sets or pairwise comparisons in multiple subsets [16-18]. Another study developed a procedure to calculate multi-set intersections by using the R software package, *SuperExactTest*, but mainly based on enumerating method [19], and the distribution theory is not established. The work also discussed the difficulty of this problem caused mainly by huge calculations, calling for a distribution theory to avoid unpractical calculations.

## 2   Definition for GHGD

In classical HGD, the parameters for the distribution include the number of the space *N*, number of the two subsets *M*, *K*, which we will denote as *M[0]*, *M[1]* for convenience (the distribution is denoted as $p(X = k \mid N, M, K) = p(X = k \mid N, M[0], M[1]) = p(X = k \mid N, M[i]_{i=0}^{1})$), here *k* is the NOE in the two subsets. However, in GHGD, suppose there are totally $T (T \geq 2)$ subsets that have *M[0]*, *M[1]*, … *M[T-1]* elements, respectively. We denote the number of elements in these subsets as an integer array $M[0, T)$ or $M[i]_{i=0}^{T-1}$. Specifically, if all of the subsets have a same number *M*, then the integer array is denoted as $M \times T$. The elements in these subsets can be overlapped for *2* to *T* times, and to more clearly elucidate this situation, we define the **level of overlap (*LO*)** for one element in the subsets: exact number of subsets that contain this elements. All of the elements that have the same *LO* is denoted as $\{LO = t\} (0 \leq t \leq T)$, and this feature is denoted as $LO = t$. In addition, we define the set containing all elements with a *LO* equal or larger than a certain number *t* as $\{LO \geq t\}$, and this feature is described as $LO \geq t$. Therefore, the parameters for the GHGD are: number of the space *N*, numbers of selected subsets $M[0, T)$, and the feature of the overlap, $LO = t$ or $LO \geq t$ (Table 1). The random variable is $X_{LO=t}$ or $X_{LO \geq t}$ (uniformly labeled as *X*), and the probability is denoted as :

$$p\left(X = k \mid N, M[i]_{i=0}^{T-1}, LO = t\right)$$

and

$$p\left(X = k \mid N, M[i]_{i=0}^{T-1}, LO \geq t\right)$$

Table 1. Definition differences between classical HGD and GHGD

| Parameters | Classical HGD | GHGD |
|---|---|---|
| Number of elements in the space | $N$ | $N$ |
| Number of selected subset | 2 | $T\ (\geq 2)$ |
| Number of elements in each subsets | $M[0], M[1]$ | $M[0], M[1], \ldots M[T-1]$ |
| Level of overlap | $LO = 2$ | $LO = t$ or $LO \geq t\ (0 \leq t \leq T)$ |
| Distribution | $p(X = k \mid N, M[0], M[1])$ $= p\left(X = k \mid N, M[i]_{i=0}^{1}\right)$ | $p\left(X = k \mid N, M[i]_{i=0}^{T-1}, LO = t\right)$ or $p\left(X = k \mid N, M[i]_{i=0}^{T-1}, LO \geq t\right)$ |

Compared to classical HGD, there is an additional parameter in GHGD: $LO = t$ or $LO \geq t$ (Table 1). To simplify, if not confused, the distributions are denoted as $p(k, LO = t)$ and $p(k, LO \geq t)$. The number of all possible subsets draws is simplified as $C(k, LO = t)$, and the mathematical expectation and variance as $E(k, LO = t)$ and $Var(k, LO = t)$. Similar simplification is used for the distributions of $p(k, LO \geq t)$.

The total numbers of all possible $T$ subsets selections from the $N$-elements space is labeled as $\left|\left\{C_N^{M[0,T)}\right\}\right|$ (the symbol $|\ |$ indicates the cardinality of a set), which can be calculated by:

$$\left|\left\{C_N^{M[0,T)}\right\}\right| = \prod_{i=0}^{T-1} C_N^{M[i]}$$

Therefore, the probability can be calculated by: $p(k, LO = t) = \dfrac{C(k, LO = t)}{\left|\left\{C_N^{M[0,T)}\right\}\right|}$, and

$p(k, LO \geq t) = \dfrac{C(k, LO \geq t)}{\left|\left\{C_N^{M[0,T)}\right\}\right|}$.

## 3 Distribution of fully overlapped elements $p(k, LO = T)$

For the simplest situation, we first considered the distribution of $p(k, LO = T)$, where $k$ represents the number of elements that are overlapped in all of the subsets. We discovered graceful formulas for essential statistics of $p(k, LO = T)$, including mathematical expectation, variance, and high order moments.

Because $C(k, LO = T) = C(k, LO \geq T)$, it was further simplified as $C(k,T)$ and can be calculated with the following recursive formula (here $M_{min}$ is the minimum of all $M[i]$ ( $0 \leq i < T$ ); $C_m^n = \dfrac{m!}{n!(m-n)!}$ is the binomial coefficient. Note: $C_m^k = 0$ when $k > m$; $C_m^0 = 1$ ):

$$C(k,T) = C_N^k \times \prod_{i=0}^{T-1} C_{N-k}^{M[i]-k} - \sum_{i=k+1}^{M_{min}} C_i^k \times C(i,T)$$

When T=2, the distribution reduced to classical HGD, confirming the correctness of the formula. The distribution is unimodal and close to symmetrical when the mean is larger enough than 0 (See Figure 1 for examples).

**Lemma 1**: $C_{N,M}(i,T) = \dfrac{N}{i} C_{N-1,M-1}(i-1,T)$. Here $C_{N,M}(k,T)$ have a distribution of $p\left(X = k \mid N, M[i]_{i=0}^{T-1}, LO = T\right)$ and $C_{N-1,M-1}(k,T)$ have a distribution of $p\left(X = k \mid N-1, (M[i]-1)_{i=0}^{T-1}, LO = T\right)$, here $(M[i]-1)_{i=0}^{T-1}$ is an integer array $M[0]-1, M[1]-1, ..., M[T-1]-1$.

Proof by mathematical induction:

First, take $i = M_{min}$, then

$$C(M_{min},T) = C_N^{M_{min}} \times \prod_{i=1}^{T} C_{N-M_{min}}^{M[i]-M_{min}} = \dfrac{N}{M_{min}} C_{(N-1)}^{(M_{min}-1)} \times \prod_{i=1}^{T} C_{(N-1)-(M_{min}-1)}^{(M[i]-1)-(M_{min}-1)} = \dfrac{N}{M_{min}} C_{N-1,M-1}(M_{min}-1,T)$$

If it is true when $i \geq k+1$, then for $i = k$ we have:

$$C(k,T) = C_N^k \times \prod_{i=0}^{T-1} C_{N-k}^{M[i]-k} - \sum_{i=k+1}^{M_{\min}} C_i^k \times C(i,T)$$

$$= \frac{N}{k} C_{N-1}^{k-1} \times \prod_{i=1}^{T} C_{N-1-(k-1)}^{M[i]-1-(k-1)} - \sum_{i=k+1}^{M_{\min}} \frac{i}{k} C_{i-1}^{k-1} \times \frac{N}{i} C_{N-1,M-1}(i-1,T)$$

$$= \frac{N}{k} \left( C_{(N-1)}^{(k-1)} \times \prod_{i=1}^{T} C_{(N-1)-(k-1)}^{(M[i]-1)-(k-1)} - \sum_{i=k+1}^{M_{\min}} C_{(i-1)}^{(k-1)} C_{N-1,M-1}(i-1,T) \right)$$

$$= \frac{N}{k} \left( C_{(N-1)}^{(k-1)} \times \prod_{i=1}^{T} C_{(N-1)-(k-1)}^{(M[i]-1)-(k-1)} - \sum_{i=(k-1)+1}^{M_{\min}-1} C_i^{(k-1)} C_{N-1,M-1}(i,T) \right)$$

$$= \frac{N}{k} C_{N-1,M-1}(k-1,T)$$

**Theorem 1**: the mathematical expectation of $p(k,T)$ is:

$$E(k,T) = \frac{\prod_{i=0}^{T-1} M[i]}{N^{T-1}}$$

Proof :

$$E(k,T) = \sum_{k=0}^{M_{\min}} k \times \frac{C_N^k \times \prod_{i=0}^{T-1} C_{N-k}^{M[i]-k} - \sum_{i=k+1}^{M_{\min}} C_i^k \times C(i,T)}{\prod_{i=0}^{T-1} C_N^{M[i]}}$$

$$= \sum_{k=0}^{M_{\min}} k \times \frac{\frac{N}{k} C_{N-1}^{k-1} \times \prod_{i=0}^{T-1} C_{N-k}^{M[i]-k} - \sum_{i=k+1}^{M_{\min}} \frac{i}{k} C_{i-1}^{k-1} \times \frac{N}{i} C_{N-1,M-1}(i-1,T)}{\prod_{i=0}^{T-1} \frac{N}{M[i]} \prod_{i=0}^{T-1} C_{N-1}^{M[i]-1}}$$

$$= \frac{\prod_{i=0}^{T-1} M[i]}{N^T} \sum_{k=0}^{M_{\min}} \frac{k \times \frac{N}{k} C_{N-1}^{k-1} \times \prod_{i=0}^{T-1} C_{N-1-(k-1)}^{M[i]-1-(k-1)} - \sum_{i=k+1}^{M_{\min}} k \times \frac{N}{k} C_{i-1}^{k-1} \times C_{N-1,M-1}(i-1,T)}{\prod_{i=0}^{T-1} C_{N-1}^{M[i]-1}}$$

$$= \frac{\prod_{i=0}^{T-1} M[i]}{N^{T-1}} \sum_{k=0}^{M_{\min}} \frac{C_{N-1}^{k-1} \times \prod_{i=0}^{T-1} C_{N-1-(k-1)}^{M[i]-1-(k-1)} - \sum_{i=k+1}^{M_{\min}} C_{i-1}^{k-1} \times C_{N-1,M-1}(i-1,T)}{\prod_{i=0}^{T-1} C_{N-1}^{M[i]-1}}$$

The sum formula is the sum of another distribution $p\left(X = k \mid N-1, (M[i]-1)_{i=0}^{T-1}, LO = T\right)$, therefore its value is 1, and the result is $\dfrac{\prod_{i=0}^{T-1} M[i]}{N^{T-1}}$.

Notably, when all M[i] is the same value M, then $E_{M[i]=M}(k,T) = \dfrac{M^T}{N^{T-1}}$.

**Theorem 2**: Any moments (*v*th) of $p(k,T)$ can be calculated by the recursion formula (considering that $E(k^0, T) = 1$):

$$E(k^v, T) = E(k, T) \times \sum_{i=0}^{v-1} C_{v-1}^i E_{N-1, M-1}(k^i, T)$$

$E_{N-1, M-1}(k^i, T)$ indicates the *i*th moments of distribution $p\left(k \mid k = |\{LO = T\}|, N-1, (M[i]-1)_{i=0}^{T-1}\right)$.

Proof:

First, when $v = 0$, $E(k^0, T) = 1$.

If it is true for *j*th (j<=v-1) moments, then:

$$E(k^v, T) = \sum_{k=0}^{M_{\min}} k^v \times \frac{C_N^k \times \prod_{i=0}^{T-1} C_{N-k}^{M[i]-k} - \sum_{i=k+1}^{M_{\min}} C_i^k \times C(i,T)}{\prod_{i=0}^{T-1} C_N^{M[i]}}$$

$$= \frac{\prod_{i=0}^{T-1} M[i]}{N^{T-1}} \times$$

$$\sum_{k=0}^{M_{\min}} ((k-1)+1)^{v-1} \times \frac{k \times \frac{N}{k} C_{N-1}^{k-1} \times \prod_{i=0}^{T-1} C_{N-1-(k-1)}^{M[i]-1-(k-1)} - \sum_{i=k+1}^{M_{\min}} k \times \frac{N}{k} C_{i-1}^{k-1} \times C_{N-1, M-1}(i-1, T)}{\prod_{i=0}^{T-1} C_{N-1}^{M[i]-1}}$$

$$= E(k, T) \sum_{k=0}^{M_{\min}} \sum_{j=0}^{v-1} C_{v-1}^j (k-1)^j \times \frac{C_{N-1}^{k-1} \times \prod_{i=0}^{T-1} C_{N-1-(k-1)}^{M[i]-1-(k-1)} - \sum_{i=k+1}^{M_{\min}} C_{i-1}^{k-1} \times C_{N-1, M-1}(i-1, T)}{\prod_{i=0}^{T-1} C_{N-1}^{M[i]-1}}$$

$$= E(k, T) \sum_{j=0}^{v-1} C_{v-1}^j \left( \sum_{k=0}^{M_{\min}} (k-1)^j \times \frac{C_{N-1}^{k-1} \times \prod_{i=0}^{T-1} C_{N-1-(k-1)}^{M[i]-1-(k-1)} - \sum_{i=k+1}^{M_{\min}} C_{i-1}^{k-1} \times C_{N-1, M-1}(i-1, T)}{\prod_{i=0}^{T-1} C_{N-1}^{M[i]-1}} \right)$$

$$= E(k, T) \sum_{j=0}^{v-1} C_{v-1}^j E_{N-1, M-1}(k^j, T)$$

**Corollary 1**: Any moments (*v*th) of $p(k,T)$ has the following formula:

$$E(k^v) = \sum_{i=1}^{v} g(v,i) \frac{\prod_{j=0}^{v-1} \prod_{i=0}^{T-1}(M[i]-j)}{\prod_{j=0}^{v-1}(N-j)^{T-1}}$$

Here $g(v,i)$ are constant coefficients which are determined by the above formulas.

**Corollary 2**: Any central moments ($v$th) about the mean can be calculated by:

$$E\left((k-E(k))^v, T\right) = \sum_{j=0}^{v} C_v^j (-1)^{v-j} E(k^j, T)(E(k,T))^{v-j}$$

**Corollary 3**: let $T=2$, then we can get any moment and central moment for classical HGD.

**Corollary 4**: let $v=2$, we can get the variance for the $p(k, LO=T)$:

$$var(k, LO=T) = E(k^2, T) - (E(k,T))^2$$

$$= E(k)\left(1 + E_{N-1, M-1}(k) - E(k)\right)$$

$$= \frac{\prod_{i=0}^{T-1} M[i]}{N^{T-1}} \times \left(1 + \frac{\prod_{i=0}^{T-1}(M[i]-1)}{(N-1)^{T-1}} - \frac{\prod_{i=0}^{T-1} M[i]}{N^{T-1}}\right)$$

Calculation of higher moments and central moments significantly characterized the intrinsic property of the distribution $p(k,T)$, and is useful in probability estimation by taking advantage of different forms of Chebyshev's inequality.

## 4  Distribution of partially overlapped elements $p(k, LO=t)$ and $p(k, LO \geq t)$ (t<T)

Calculation of $C(k, LO=t)$ and $C(k, LO \geq t)$ is much more complicated than that of $C(k,T)$. Although there is no simple formula, it can be calculated by recursive algorithms implemented with program code (the code written by java is included in supplementary file 1, ProgramMethods.java). The essential statistics of the distribution including mathematical expectation and variance also have much more complicated forms (Table 2), which can be calculated by formulas when the set number is no more than 7 ($T \leq 7$), while for larger $T$ the corresponding formulas can be obtained with the help of more powerful computers.

### 4.1  Programmable algorithm for $p(k, LO=t)$ and $p(k, LO \geq t)$ calculations

The key concept for calculation of $C(k, LO=t)$ and $C(k, LO \geq t)$ is as follows: select one subset $M[i-1]$ each time, and store the state of $|\{LO=t\}|$ ($0 \leq t \leq i$) each time in array

$R[i-1][t]$ (here $0 \leq t \leq i$). That is, the N-elements set was segmented into $i+1$ subsets each has different LO ($0 \leq LO \leq i$). Do next subset draw of selecting M[i] elements from the N-elements set, supposing selecting $Km[j]$ ($0 \leq j \leq i$ and $\sum_{j=0}^{i} Km[j] = M[i]$) elements in each segmented subset $\{LO = t\}$ ($0 \leq t \leq i$), then the state of $R[i][t]$ is updated by:

$$\begin{cases} R[i][0] = R[i-1][0] - Km[0] \\ R[i][t] = R[i-1][t] - Km[t] + Km[t-1] \quad 0 < t < i \\ R[i][i+1] = Km[i] \end{cases}$$

At the end of recursive calculation, the total number for each segmented T+1 subsets with different LO ($0 \leq LO \leq T$) was stored in $R[T][T+1]$, which can be used to calculate $C(k, LO = t)$ and $C(k, LO \geq t)$ ($0 \leq t \leq T$). From the calculated results, the distribution is unimodal (Figure 1), and the distribution feature depends on the relative values of N and M.

### 4.2 Mathematical expectation and variance of the distributions $p(k, LO = t)$ and $p(k, LO \geq t)$ ($0 \leq t \leq T$)

By performing huge amount of iterative and extensive induction and analysis, we got explicit formulas for the mathematical expectation and variance for most $p(k, LO = t)$ and $p(k, LO \geq t)$ ($0 \leq t \leq T$). Briefly, serial values of $N$, $M[]$ and $T$ were taken to calculate the $p(k, LO = t)$, $E(k, LO = t)$ and $Var(k, LO = t)$ by the above programmable algorithm. $E(k, LO = t)$ and $Var(k, LO = t)$ were found to be algebraic expressions of the polynomials of $N$, $M[]$ and $T$. Then the coefficients were calculated by solving a great amount of interrelated linear equations. The most important results were presented below for convenient application in practice.

#### 4.2.1 Mathematical expectations for $p(k, LO = t)$ and $p(k, LO \geq t)$ ($0 \leq t \leq T$)

**Result 1**: The formula for $E(k, LO = t)$ and $E(k, LO \geq t)$ are:

$$E(k, LO = t) = \sum_{l=0}^{T-t} C_{t+l}^{l} \times (-1)^{l} \times \frac{S(M^{t+l})}{N^{t+l-1}}$$

And

$$E(k, LO \geq t) = \sum_{l=0}^{T-t} C_{t+l}^{0 \to \pm l} \times \frac{S(M^{t+l})}{N^{t+l-1}}$$

$S(M^z)$ is the elementary symmetric polynomial of M[0, T) which is defined as the sum of all products of $z$ distinct variables in M[0, T):

$$S(M^z) = \sum_{i_{j_1} \neq i_{j_2}(j_1 \neq j_2)} \prod_{j=0}^{z-1} M[i_j]$$

And

$$C_n^{m1 \to \pm m2} = \sum_{m=m1}^{m2} (-1)^m \times C_n^m$$

These formulas were validated by the above programmable algorithm and by exhaustive enumeration algorithm.

### 4.2.2 Variances for $p(k, LO = t)$ ($0 \leq t \leq T$)

The formulas for the variance of these distributions are much more complex. We will give the formula for $E(k^2, LO = t) - E(k, LO = t)$, then $Var(k, LO = t)$ can be calculated by

$Var(k, LO = t) = E(k^2, LO = t) - (E(k, LO = t))^2$.

**Result 2**: For $var(k, LO = t)$ $t \leq T$, we first got the formula for the situation that all $M[i] = M$ ($M \times T$) ($d = T - t$):

$$E_{M \times T}(k^2, LO = t) - E_{M \times T}(k, LO = t) = \sum_{x=0}^{d} \sum_{o=0}^{d} \sum_{l=o}^{t+o} \frac{(-1)^{t+l+x} C_d^o C_{d+l-o}^d C_{t+o}^l C_{t+x}^x C_T^{t+x} M^{l+t+x}}{(N-1)^{t+o-1} N^{t+x-1}}$$

And for the general situation of $M[i]$, we have ($d = T - t$):

$$E(k^2, LO = t) - E(k, LO = t) = \sum_{x=0}^{d} \sum_{o=0}^{d} \sum_{l=o}^{t+o} \frac{(-1)^{t+l+x} C_d^o C_{d+l-o}^d C_{t+o}^l C_{t+x}^x C_T^{t+x} S^2(M^{t+l+x})}{(N-1)^{t+o-1} N^{t+x-1}}$$

Here $S^2(M^z)$ is 2-order symmetric polynomials of $M[0,T)$, that is the sum of distinct products of powers (<=2) of z distinct variables in $M[0,T)$. The general calculation of

$S^2(M^z)$ is not known yet. Nevertheless, all of the $S^2(M^z)$ is theoretically calculable and the coefficients for all $T \leq 7$ have been deduced for convenient use (supplementary file 2, MCoesT1To7-2.txt).

### 4.2.3 Variances for $p(k, LO \geq t)$ ($0 \leq t \leq T$)

The formula of $Var(k, \geq t)$ work out to rather complicated expressions. We have the following result:

$$E(k^2, LO \geq t) = \left(1 + \frac{F(M)}{X(M) \times (N-1)^{T-1}}\right) \times E(k, LO \geq t)$$

Where $F(M)$ and $X(M)$ are polynomial functions for $N$, $M[]$, $T$ and $t$, and can be calculated as described below.

Table 2. Comparison of essential statistics for HGD and GHGD. The similarities between the formulas of mathematical expectation and variance for HGD and GHGD explains why GHGD is a general form for HGD.

| Statistics | HGD* | GHGD ($LO = T$) |
|---|---|---|
| Random variable | NOE that are overlapped in the two subsets | NOE that are overlapped exactly $t$ or at least $t$ times in the $T$ subsets |
| Mathematical expectation | $\frac{M[0]M[1]}{N}$ | $\frac{\prod_{i=0}^{T-1} M[i]}{N^{T-1}}$ |
| Variance | $\frac{M[0]M[1]}{N}\left(1 + \frac{(M[0]-1)(M[1]-1)}{N-1} - \frac{M[0]M[1]}{N}\right)$ # | $\frac{\prod_{i=0}^{T-1} M[i]}{N^{T-1}} \times \left(1 + \frac{\prod_{i=0}^{T-1}(M[i]-1)}{(N-1)^{T-1}} - \frac{\prod_{i=0}^{T-1} M[i]}{N^{T-1}}\right)$ |
| High central moments ($v$th) | | $\sum_{j=0}^{v} C_v^j (-1)^{v-j} E(k^j, T)(E(k,T))^{v-j}$ $E(k^v, T) = E(k,T) \times \sum_{i=0}^{v-1} C_{v-1}^i E_{N-1, M-1}(k^i, T)$ |

Table 2. Continued. Essential statistics for GHGD ($LO = t$ and $LO \geq t, 1 \leq t \leq T$)

| Statistics | GHGD ($LO=t$)($1 \leq t \leq T$) | GHGD ($LO \geq t$)($1 \leq t \leq T$) |
|---|---|---|
| Mathematical expectation | $\sum_{l=0}^{T-t} C_{t+l}^l \times (-1)^l \times \frac{S(M^{t+l})}{N^{t+l-1}}$ ## | $\sum_{l=0}^{T-t} C_{t+l}^{0 \to \pm l} \times \frac{S(M^{t+l})}{N^{t+l-1}}$ † |
| $E(k^2) - E(k)$ ** | $\sum_{x=0}^{d} \sum_{o=0}^{d} \sum_{l=o}^{t+o} \frac{(-1)^{t+l+x} C_d^o C_{d+l-o}^d C_{t+o}^l C_{t+x}^x C_T^{t+x} S^2(M^{t+l+x})}{(N-1)^{t+o-1} N^{t+x-1}}$ ‡ | $\left(1 + \frac{F(M)}{X(M) \times (N-1)^{T-1}}\right) \times E(k, LO \geq t)$ <br> $F(M) = \sum_{x=0}^{T+d}(-1)^{T+d-x}\left[\sum_{y=0}^{2d} q(x,y)(N-1)^y\right] S'^2(M^{x+t})$ Δ <br> $X(M) = \sum_{x=0}^{d}(-1)^x C_{t-1+x}^{t-1} C_T^{d-x} N^{d-x} \frac{S(M^{x+t})}{C_T^{x+t}}$ ## |

Note: * Select M[1] balls from a finite population space of N balls containing M[0] red and (N-M[0]) green balls;

\# the variance for classical HGD is $\frac{M[0]M[1]}{N}\left(\frac{(N-M[0])(N-M[1])}{N(N-1)}\right) =$

$$\frac{M[0]M[1]}{N}\left(1+\frac{(M[0]-1)(M[1]-1)}{N-1}-\frac{M[0]M[1]}{N}\right); \text{\#\#} \quad S(M^z) = \sum_{i_{j_1} \neq i_{j_2}(j_1 \neq j_2)} \prod_{j=0}^{z-1} M[i_j]; \dagger$$

$C_n^{m1 \to \pm m2} = \sum_{m=m1}^{m2}(-1)^m \times C_n^m$; ‡ $S^2(M^z)$ is 2-order symmetric polynomials of $M[0,T]$, which have been deduced for all $T \leq 7$ (supplementary file 2, MCoesT1To7-2.txt); △ $S'^2(M^z)$ is similar to $S^2(M^z)$, but has different coefficients. Their values are also deduced for all $T \leq 7$ (supplementary file 2, MCoesTmoreTo7.txt); the values for $q(x,y)$ for all $T \leq 7$ is also deduced and available (supplementary file 2, qxy.txt). **Variance is calculated by: $Var(k) = E(k^2) - (E(k))^2$ for both $LO = t$ and $LO \geq t$.

**Result 3**: First, for situation that all $M[i] = M$ ($M \times T$) ($d = T - t$):

$$X_{M \times T}(M) = \sum_{x=0}^{d}(-1)^x C_{t-1+x}^{t-1} C_T^{d-x} M^x N^{d-x}$$

And

$$F_{M \times T}(M) = \sum_{x=0}^{T+d}(-1)^{T+d-x}\left(\sum_{y=0}^{2d} q(x,y)(N-1)^y\right) M^x$$

$q(x,y)$ is a function for $T$, $t$, $x$ and $y$. The general formula for $q(x,y)$ is not known, while it can be theoretically deduced and the values for all $T \leq 7$ have been obtained for convenient use (supplementary file 2, qxy.txt).

**Result 4**: For the general situation of $Var(k, LO \geq t)$, we have:

$$X(M) = \sum_{x=0}^{d}(-1)^x C_{t-1+x}^{t-1} C_T^{d-x} N^{d-x} \frac{S(M^{x+t})}{C_T^{x+t}}$$

$$F(M) = \sum_{x=0}^{T+d}(-1)^{T+d-x}\left(\sum_{y=0}^{2d} q(x,y)(N-1)^y\right) S'^2(M^{x+t})$$

Here $S'^2(M^z)$ is similar to $S^2(M^z)$, but has different coefficients. Their values are also theoretically calculable and the results for $T \leq 7$ have been deduced for practice use (supplementary file 2, MCoesTmoreTo7.txt).

## 5 Application

### 5.1 Statistical theory for Application

The mathematical expectation and variance can be feasibly calculated in most situations with the above result. Although the mode for the GHGD is not known yet, as an empirical result after checking a great amount of distributions with various *N* and *M*[0, *T*) values, the difference between mean and mode is very small and are normally less than 1 (Figure 1). Statistical inference can be made with these results and information.

Assume there is a "*theoretical set consists of truly targeted elements*" (TSTE) in the universal space. In experimental results, suppose we got several sample subsets (SS) from a space having *N* elements in experiments and each subset has $M_{sample}[i]_{i=0}^{T-1}$ elements. Each subset may partially target the TSTE, and elements with higher LO are more likely to target the TSTE (Figure 2). The elements that are overlapped exactly *t* times or at least *t* times in the SSs are denoted as "*overlapped elements in the sample sets with LO = t or LO ≥ t*" (OESS$_{LO=t}$ or OESS$_{LO≥t}$), and the number of these elements are denoted as NOESS$_{LO=t}$ and NOESS$_{LO≥t}$. To simplify, if not confused, "OESS$_{LO=t}$ or OESS$_{LO≥t}$" are collectively referred to as "OESS", and "NOESS$_{LO=t}$ or NOESS$_{LO≥t}$" as "NOESS".

The central questions are: are these OESS target the TSTE with statistical significance? And if yes, how much proportion of the OESS target the TSTE with statistical significance? By performing statistical inference for serial *t* values ($1 \leq t \leq T$), we can tell for each *t* value how much portion of the OESS can be used as valuable elements, especially in biological research to tell which genes can be used for further study.

Because the TSTE is a theoretical concept and normally not known, for statistical inference, the probability of event $p(X \geq NOESS_{LO=t} \mid N, M[i]_{i=0}^{T-1}, LO = t)$ or $p(X \geq NOESS_{LO≥t} \mid N, M[i]_{i=0}^{T-1}, LO \geq t)$ is used to estimate **the probability that "the OESS** hit the TSTE". Similarly, "$p\left(X \geq NOESS_{LO=t} \mid N, M[i]_{i=0}^{T-1}, LO=t\right)$ or $p\left(X \geq NOESS_{LO \geq t} \mid N, M[i]_{i=0}^{T-1}, LO \geq t\right)$" were simplified as "$p(X \geq NOESS)$".
Therefore, if $p(X \geq NOESS)$ is lower than a statistical level α (normally 0.05), then we can tell that the OESS hit the TSTE significantly. For quantitative analysis of how much proportion of the OESS hit the TSTE with statistical significance, we need to further analyze the probability of event $p(X > Z) < \alpha$). Let $Z_{min}$ be the smallest value of $Z$ that satisfy the condition $p(X > Z) < \alpha$. Then the number and percent of elements in the OESS that hit the TSTE (defined as *statistical hit number, SHN, and statistical hit rate,* SHR, respectively) is not smaller than $NOESS - Z_{min}$ and $\frac{NOESS - Z_{min}}{NOESS}$, respectively. Notably, if $Z_{min} <$ NOESS, then SHR>0 and OESSs hit the TSTE significantly; if $Z_{min}=0$ then SHR=100% and all of OESS hit the TSTE significantly. These two kinds of probability were denoted as $p_{SHR>0}$ and $p_{SHR=100\%}$, respectively.

The above probabilities can be estimated by using different kinds of Chebyshev's inequality ($\mu = E(k)$, and $\sigma^2$ is the variance) for unimodal distribution(Savage, 1961):

$$P(|X - \mu| \geq \lambda) \leq \frac{\sigma^2}{\lambda^2}$$

Or (here $s = |\mu - \mu_0|$, and $\mu_0$ is the mode of the distribution)

$$P(|X - \mu| \geq \lambda) \leq \frac{4(\sigma^2 + s^2)}{9(\lambda - s)^2} \quad if\ \lambda > s$$

**5.2 Example in practice**

We developed a Java software that can be easily used to perform statistical inference and to investigate the characteristics of the distributions (supplementary file 3). As an example, to infer potential interacting genes of key gene OLIG2 in glioblastoma, the most lethal primary brain cancer in adult, we got OLIG2 interacting genes from 4 independent datasets by performing ARACNe algorithm (Mao et al., 2013; Margolin et al., 2006). As a result, we got four lists of gene subsets (SSs), containing 157, 110, 87, 110 genes, respectively. Considering that there are totally about 19815 genes in the human genome used in the gene chip, we have $N$=19815, and $M$={127, 110, 87, 110}. For all of these genes, there are 14, 25, 44 and 245 genes that are overlapped 4, 3, 2, and 1 times, that is $NOESS_{LO=4}=14$, $NOESS_{LO=3}=25$, $NOESS_{LO=2}=44$, $NOESS_{LO=1}=245$, and $NOESS_{LO\geq t}=14$, $NOESS_{LO\geq 3}=39$, $NOESS_{LO\geq 2}=83$, $NOESS_{LO\geq 1}=328$ (Figure 2).

Table 3. Distribution analysis for N=19815, M={139, 109, 372, 249}

| Distribution with different $LO$ | Mathematical expectation | Variance | NOESS | $p(X\geq NOESS)$ ($p_{SHR>0}$) | $p(X\geq 1)$ ($p_{SHR=100\%}$) | 95% credibility interval ($\approx$) | SHN and SHR (p<0.05) |
|---|---|---|---|---|---|---|---|
| $p(k,t=4)$ | 0.000017 | 0.000017 | 14 | < 0.000000 | <0.000017 | [0, 0.018] | 14 (100%) |
| $p(k,t=3)$ | 0.012717 | 0.012714 | 25 | < 0.000020 | <0.013043 | [0, 0.516] | 25 (100%) |
| $p(k,t=2)$ | 3.505980 | 3.393887 | 44 | < 0.001252 | >-0.459566 | [0, 11.74] | 33 (75%) |
| $p(k,t=1)$ | 426.949821 | 13.682894 | 145 | < 0.000083 | >-0.999964 | [410.41, 443.49] | NULL |
| $p(k,t\geq 4)$ | 0.000017 | 0.000017 | 14 | < 0.000000 | <0.000017 | [0, 0.018] | 14 (100%) |
| $p(k,t\geq 3)$ | 0.012734 | 0.012731 | 39 | < 0.000008 | <0.013062 | [0, 0.517] | 39 (100%) |
| $p(k,t\geq 2)$ | 3.518714 | 3.405393 | 83 | < 0.000318 | >-0.463204 | [0, 11.77] | 72 (86.75%) |
| $p(k,t\geq 1)$ | 430.468535 | 3.442445 | 328 | < 0.000192 | >-0.999989 | [422.17, 438.76] | NULL |

First, by applying the Chebyshev's inequality, all of the OESS with LO ≥ 2 hit the TSTE with statistical significance (Table 3; $p_{SHR>0}=p(X\geq NOESS)$). Next, we investigated whether all of the OESS with certain LO are significantly hit the TSTE by calculating the value of $p_{SHR=100\%}=p(X\geq 1)$. In practice, we are more interested in elements that are overlapped **at least** $t$ times, and we will focus our discussion on $p(k,LO\geq t)$. As a result, $p(|X|\geq 1,LO\geq 4)\leq 0.000017$, and $p(|X|\geq 1,LO\geq 3)\leq 0.013$. If the 4 sets were randomly selected from the genome, the probability of at most one gene ($|X|\geq 1$) overlapped in at least 3 sets are $\leq 0.013 < 0.05$. Therefore, all of the genes (totally 14+25=39) with $LO\geq 3$ significantly hit the TSTE. However, for the situation of $p(k, LO \geq 2)$, the 95% credibility interval around the mean is about [0, 11.74], implying that, if randomly selected, there would be at most 11 genes with $LO \geq 2$. Now we have 83 genes with $LO \geq 2$ in OESS$_{LO \geq 2}$, indicating a portion of these genes (at least totally $83 - 11 = 72, 86.75\%$) hit the TSTE at a statistical level of $p = 0.05$ (hence SHR= 86.75%).

As conclusions for this example (Table 3): 1. All of the OESS with $LO \geq 2$, 3 or 4 hit the TSTE with statistical significance ($p < 0.05$); 2. All of the fully overlapped genes (*LO* = 4) are significantly overlapped and hit the TSTE with $p \leq 0.000017$; 3. All of the genes overlapped in at least 3 subsets ($LO \geq 3$) are significantly overlapped and hit the TSTE with $p \leq 0.013$; 4. A portion of the genes overlapped in at least 2 sets (86.75%, $LO \geq 2$) hit the TSTE with $p < 0.05$.

## 6 Conclusion

By generalizing the classical HGD, we established a statistical theory for the distribution of NOE in multiple subsets selected from a finite population. According to the results, essential statistics of the distribution including mathematical expectation and variance can be easily calculated by explicit formulas and used to perform statistical inference. For the fully overlapped elements ($LO = T$), all essential statistics were given including mathematical expectations, variance, and high order moments. For partially overlapped elements ($LO = t$ or $LO \geq t$, $1 \leq t < T$), the mathematical expectations are known, while variances can be calculated for $T \leq 7$ with explicit formulas (and for larger *T* it can be deduced with the help of more powerful computers). Although the concept of GHGD is inspired from biological research, it would be useful in a wide range of situations involving subsets selection including physics and mathematical research.

**ACKNOWLEDGMENTS**

Funding and Conflict of interest: This study was partially supported by National Natural Science Foundation of China (81502143).

**Supplementary materials:**

Supplementary Table S1.

Supplementary File 1-3.

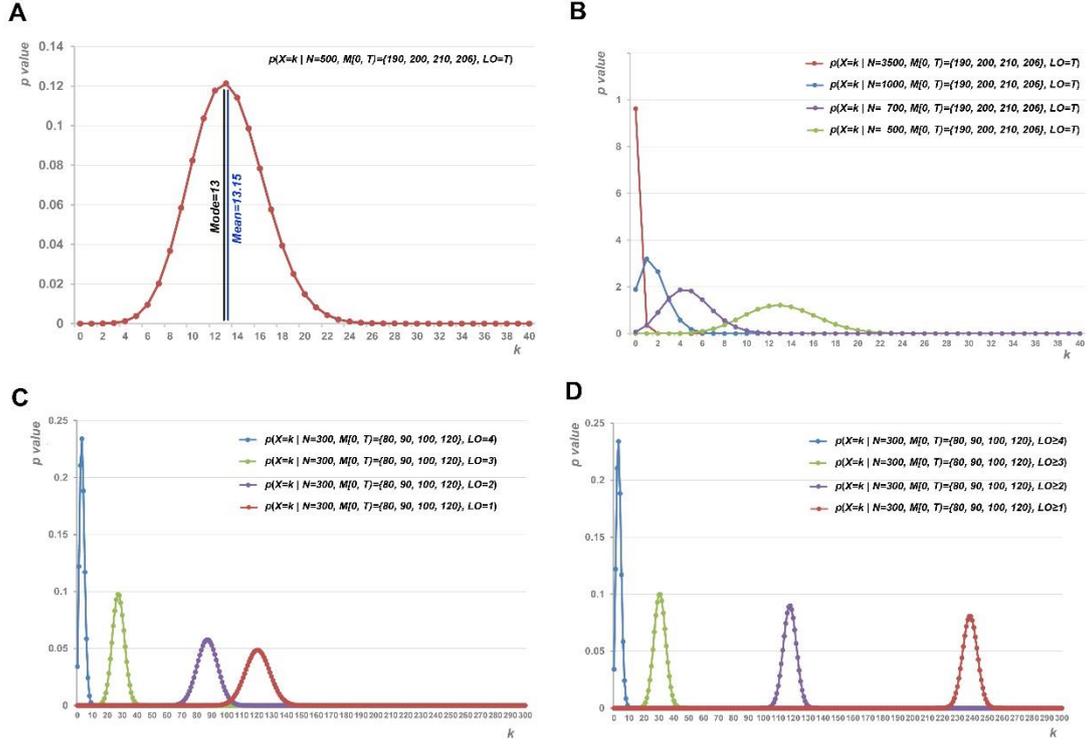

Figure 1. Example distributions of $p(k,T)$, $p\left(X = k \mid N, M[i]_{i=0}^{T-1}, LO = t\right)$ and $p\left(X = k \mid N, M[i]_{i=0}^{T-1}, LO \geq t\right)$, illustrating different shapes of the distribution curve. (**A**). A typical GHGD is unimodal. Note the distribution has a mode very close to the mean (mode = 13, mean = 13.15). (**B**). The influence of relative *N* and *M* values on the GHGD curves. (**C-D**). All of the distribution curves (for all $LO = t$ (C) and $LO \geq t$ (D), $1 \leq t \leq T$) for the same *N* and *M* value.

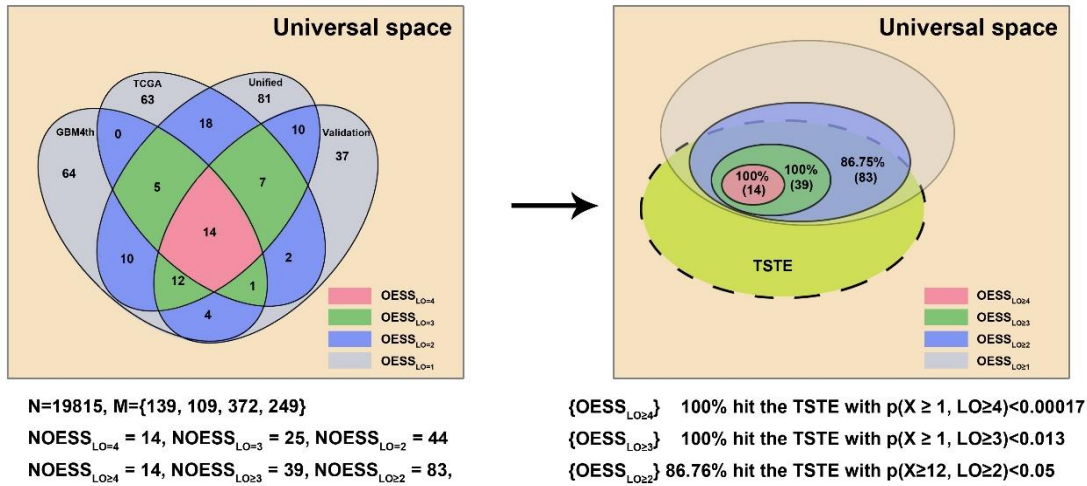

Fig. 2 Vein diagram and statistical inference for overlapped OLIG2-interacting-genes inferred from 4 different datasets. $OESS_{LO=t}$ or $OESS_{LO\geq t}$: *overlapped elements in the sample sets with $LO = t$ or $LO \geq t$* ; $NOESS_{LO=t}$ or $NOESS_{LO\geq t}$: the number of $OESS_{LO=t}$ or $OESS_{LO\geq t}$; TSTE: *theoretical set consists of truly targeted elements* in the universal space.

**Supplementary materials:**

Supplementary Table S1: List of a small portion of recently published papers utilizing the Venn Diagrams to get overlapped elements.

Supplementary File 1: Java code for the Programmable algorithms for calculation of $C(k, LO=t)$ and $C(k, LO \geq t)$.

Supplementary File 2: Lists of coefficients of $S^2(M^z)$, $S'^2(M^z)$ and $q(x,y)$ for the calculation of $Var(k, LO=t)$ and $Var(k, LO \geq t)$ (related to result 2 and 3 in part "*Variances for $p(k, LO=t)$ ( $0 \leq t \leq T$ )*" and "*Variances for $p(k, LO \geq t)$ ( $0 \leq t \leq T$ )*").

Supplementary File 3: Runnable java software (jar file) for investigation and statistical inference for $p(k, LO=t)$ and $p(k, LO \geq t)$ distributions, and demo data files (interacting genes of OLIG2 inferred from 4 datasets).